\documentclass{amsart}

\usepackage{amsmath,amsthm,amsfonts,amssymb}

\usepackage[backref=page]{hyperref}
\hypersetup{colorlinks=true,citecolor=blue,linkcolor=blue,urlcolor=blue,
pdfstartview=FitH, pdfauthor=Valentin Blomer and Gergely
Harcos,pdftitle=The spectral decomposition of shifted convolution
sums}

\theoremstyle{plain}
\newtheorem{theorem}{Theorem}

\theoremstyle{definition}
\newtheorem{remark}{Remark}

\renewcommand{\geq}{\geqslant}
\renewcommand{\leq}{\leqslant}

\DeclareMathOperator{\PGL}{PGL} \DeclareMathOperator{\SL}{SL}
\DeclareMathOperator{\ord}{ord} \DeclareMathOperator{\sgn}{sgn}
\DeclareMathOperator{\vol}{vol} \DeclareMathOperator{\Ad}{Ad}
\DeclareMathOperator*{\res}{res}

\newcommand{\eps}{\varepsilon}
\newcommand{\GG}{\Gamma}
\newcommand{\RR}{\mathbb{R}}
\newcommand{\CC}{\mathbb{C}}
\newcommand{\ZZ}{\mathbb{Z}}
\newcommand{\DD}{\mathcal{D}}
\newcommand{\HH}{\mathcal{H}}
\newcommand{\KK}{\mathcal{K}}
\newcommand{\fg}{\mathfrak{g}}

\begin{document}

\author{Valentin Blomer}
\address{Department of Mathematics, University of Toronto, 40 St. George Street, Toronto, Ontario, Canada, M5S 2E4} \email{vblomer@math.toronto.edu}

\author{Gergely Harcos}
\address{Alfr\'ed R\'enyi Institute of Mathematics, Hungarian Academy of Sciences, POB 127, Budapest H-1364, Hungary} \email{gharcos@renyi.hu}
\title{The spectral decomposition of shifted convolution sums}

\thanks{The first author was in part supported by an NSERC research grant 311664-05.
The second author was supported by NSF grant DMS-0503804 and
European Community grants MTKD-CT-2004-002988 and
MEIF-CT-2006-040371 under the Sixth Framework Programme.}

\keywords{shifted convolution sums, spectral decomposition, Hecke eigenvalues, non-holomorphic cusp forms, Kirillov model}

\begin{abstract} Let $\pi_1$, $\pi_2$ be cuspidal automorphic representations
of $\PGL_2(\RR)$ of conductor $1$ and Hecke eigenvalues
$\lambda_{\pi_{1, 2}}(n)$, and let $h>0$ be an integer. For any
smooth compactly supported weight functions $W_{1,
2}:\RR^\times\to\CC$ and any $Y>0$ a spectral decomposition of the
shifted convolution sum
\begin{displaymath}
  \sum_{m\pm n=h}\frac{\lambda_{\pi_1}(|m|)\lambda_{\pi_2}(|n|)}{\sqrt{|mn|}}
W_1\left(\frac{m}{Y}\right)W_2\left(\frac{n}{Y}\right)
\end{displaymath}
is obtained. As an application, a spectral decomposition of the Dirichlet series
\begin{displaymath}
  \sum_{\substack{m,n\geq 1\\ m-n=h}}
\frac{\lambda_{\pi_1}(m)\lambda_{\pi_2}(n)}{(m+n)^{s}} \left(\frac{\sqrt{mn}}{m+n}\right)^{100}
\end{displaymath}
is proved for $\Re s > 1/2$ with polynomial growth on vertical lines
in the $s$ aspect and uniformity in the $h$ aspect.
\end{abstract}

\subjclass[2000]{Primary 11F30, 11F70, 11F72; Secondary 11F12,
11M41}

\maketitle
\tableofcontents

\section{Introduction}

Let $G=\PGL_2(\RR)$ and $\GG=\PGL_2(\ZZ)$. There is a spectral
decomposition
\begin{equation}\label{spectraldec}L^2(\GG\backslash G)=\int_\tau
V_\tau\,d\tau,\end{equation} where $(\tau,V_\tau)$ are irreducible
automorphic representations of $\GG\backslash G$ (including the
trivial representation) and $d\tau$ is the spectral measure defined
as follows: The trivial representation $(\tau_0,\CC)$ has spectral
measure $1$. Each nontrivial representation $(\tau,V_\tau)$ is
generated by a modular form on $\HH$ with respect to the full
modular group $\SL_2(\ZZ)$, hence it corresponds to some Laplacian
eigenvalue
\begin{displaymath}
  \lambda_\tau=1/4-\nu_\tau^2\in\RR \quad\text{with}\quad\Re\nu_t\geq 0,\ \Im\nu_\tau\geq 0.
\end{displaymath}
We shall use the notation
\begin{displaymath}
 \tilde{ \lambda}_{\tau} := 1 + |\lambda_{\tau}|.
\end{displaymath}
A representation $(\tau,V_\tau)$ generated by a holomorphic or Maass
cusp form has spectral measure $1$. We can assume that the
underlying cusp form is a Hecke eigenform, and we denote by
$\lambda_\tau(n)$ its $n$-th Hecke eigenvalue. A representation
$(\tau,V_\tau)$ generated by an Eisenstein series with
$\Re\nu_\tau=0$ and $\Im\nu_\tau>0$ has spectral measure
$d\nu_\tau/(2\pi i)$, and we denote by $\lambda_\tau(n)$ the divisor
sum $\sum_{ab=n}(a/b)^{\nu_\tau}$.

For a function $W \in \mathcal{C}^d(\RR^\times)$ we denote
\begin{equation}\label{akshaynorm}\|W\|_{A^d}:=\sum_{j=0}^{d}
\left(\int_{\RR^\times}(|u|+|u|^{-1})^{d}\left|\frac{
d^jW}{du^j}\right|^2d^\times u\right)^{1/2},\end{equation} provided
the integral is finite. Here $d^\times u:= du/|u|$ is the Haar
measure on $\RR^\times$.\\

In this paper we obtain a spectral decomposition for shifted
convolution sums of Hecke eigenvalues of two arbitrary cusp forms,
as well as a spectral decomposition of the corresponding Dirichlet
series with polynomial growth estimates on vertical lines and
uniform dependence with respect to the shift parameter.

\begin{theorem}\label{theorem1}
Let $\pi_1$ and $\pi_2$ be arbitrary cuspidal automorphic
representations of $\GG\backslash G$. Let $a, b, c \geq 0$ be
arbitrary integers and let $W_{1,2}:\RR^\times\to\CC$ be arbitrary
functions such that $\|W_{1, 2}\|_{A^d}$ exist for $d = 18
+2a+2b+4c$. Then there exist functions $W_{\tau}:(0,\infty)\to\CC$
depending only on $\pi_{1,2}$, $W_{1,2}$, and $\tau$ such that the
following two properties hold:

$\bullet$ If $h>0$ is an arbitrary integer and $Y>0$ is arbitrary
then one has the decomposition over the full spectrum (excluding the
trivial representation)
\begin{equation}\label{spectraldecomposition}\sum_{m+n=h}\frac{\lambda_{\pi_1}(|m|)\lambda_{\pi_2}(|n|)}{\sqrt{|mn|}}
W_1\left(\frac{m}{Y}\right)W_2\left(\frac{n}{Y}\right)=
\int_{\tau\neq\tau_0}\frac{\lambda_\tau(h)}{\sqrt{h}}W_{\tau}\left(\frac{h}{Y}\right)d\tau.\end{equation}

$\bullet$ If $0<\eps<1/4$ is arbitrary then for $y > 0$ one has the
uniform bound
\begin{equation}\label{ineq9}\int_{\tau\neq\tau_0}\tilde{\lambda}_\tau^{c}
\left|\left(y\frac{d}{dy}\right)^{a}W_{\tau}(y)\right|\,d\tau
\ll_{\eps,a,b,c}
C_{a,b,c}\,\min(y^{1/2-\eps},y^{1/2-b-\eps}),\end{equation} where
\[C_{a,b,c}:=(\tilde{\lambda}_{\pi_1}+\tilde{\lambda}_{\pi_2})^{11+a+b+2c}
\sum_{d_1+d_2=18+2a+2b+4c}\|W_1\|_{A^{d_1}}\,\|W_2\|_{A^{d_2}}.\]
\end{theorem}

\begin{remark} The condition on $W_{1,2}$ is satisfied as
long as these functions decay rapidly near $0$ and $\pm \infty$. One
can view $W$ as a function on $(0, \infty)$ times (a subset of) the
unitary dual of $\GG\backslash G$. Implicit in
\eqref{spectraldecomposition} is the fact that this function is
integrable in the second variable with respect to the spectral
measure. By explicating $W$ stronger regularity properties would
follow. Finally, the right-hand side of
\eqref{spectraldecomposition} captures also the summation condition
$m -n =h$ if the support of $W_2$ is on the negative axis.
\end{remark}

By a result of Kim--Sarnak \cite[Appendix~2]{KS} we have
\begin{displaymath}
  |\lambda_{\tau}(h)| \leq d(h) h^{\theta}
\end{displaymath}
for $\theta = 7/64$ and any cuspidal automorphic representation
$\tau$ of $\GG\backslash G$, where $d(n)$ is the number of divisors
of $n$. Thus for two smooth compactly supported functions $W_{1,
2}$, the bound \eqref{ineq9} gives immediately
\begin{equation}\label{shiftedbound}
  \sum_{m \pm n = h} \lambda_{\pi_1}(|m|) \lambda_{\pi_2}(|n|) \,W_1\left(\frac{m}{Y}\right)W_2\left(\frac{n}{Y}\right)
  \ll_{\eps,\pi_1, \pi_2,W_1, W_2} h^{\theta}Y^{1/2} (hY)^{\eps}
\end{equation}
for any $\eps > 0$, uniformly in $h>0$, cf.\ e.g.\ \cite{Bl}. The
novelty of Theorem~\ref{theorem1} is to obtain an exact spectral
decomposition of the left side of \eqref{spectraldecomposition}
rather than an upper bound. Thus Theorem~\ref{theorem1} develops its
full strength when the left side of \eqref{spectraldecomposition} is
averaged over $h$, as is necessary, for example, to prove subconvex
estimates for certain families $L$-functions. If $\pi$ is a cuspidal
automorphic representation of $G$ of arbitrary conductor and $\chi$
is a primitive Dirichlet character of conductor $q$, then without
much effort we can deduce the Burgess-type bound
\begin{equation}\label{burgessbound}
  L(1/2,\pi \otimes \chi) \ll_{\pi, \eps}  q^{\frac{3}{8} + \frac{\theta}{4}+\eps}
\end{equation}
by combining a slight generalization of Theorem~\ref{theorem1} with
amplification and a large sieve inequality (cf.\
Remark~\ref{remark2} below). In fact, a similar result can be
derived more generally over totally real number fields, where
classical methods are much harder to implement. We postpone this
discussion to \cite{BH} and give here another application of
Theorem~\ref{theorem1} that was our initial motivation.

\begin{theorem}\label{theorem2} Let $\pi_1$ and $\pi_2$ be arbitrary
cuspidal automorphic representations of $\GG\backslash G$, and let
$c,k\geq 0$ be arbitrary integers satisfying $k>60+12c$. There are
holomorphic functions $F_{k,\tau}:\{s:\frac{1}{2}<\Re
s<\frac{3}{2}\}\to\CC$ depending only on $\pi_{1,2}$, $c$, $k$, and
$\tau$ such that the following two properties hold:

$\bullet$ If $h>0$ is an arbitrary integer then one has the
decomposition over the full spectrum (excluding the trivial
representation)
\begin{equation}\label{spectraldecomposition2}\sum_{\substack{m,n\geq 1\\ m-n=h}}
\frac{\lambda_{\pi_1}(m)\lambda_{\pi_2}(n)\,(mn)^\frac{k-1}{2}}{(m+n)^{s+k-1}}
=h^{\frac{1}{2}-s}\int_{\tau\neq\tau_0}\lambda_\tau(h)\,F_{k,\tau}(s)\,d\tau,\qquad
\Re s>1.\end{equation}

$\bullet$ One has the uniform bound
\begin{equation}\label{ineq15}\int_{\tau\neq\tau_0}\tilde{\lambda}_\tau^{c}\,|F_{k,\tau}(s)|\,d\tau\ll_{\eps,k}
(\tilde{\lambda}_{\pi_1}+\tilde{\lambda}_{\pi_2})^{12+4c}\,|s|^{22+4c},\quad
\frac{1}{2}+\eps<\Re s<\frac{3}{2}.\end{equation}
\end{theorem}

\begin{remark}\label{remark2}
Theorems~\ref{theorem1} and \ref{theorem2} extend to arbitrary level
and to the more general additive constraints $\ell_1 m\pm\ell_2 n=h$
in a straightforward fashion with good control in the $\ell_{1,2}$
parameters. Here small technical complications arise from the
possible presence of complementary series representations $\tau$ and
the presence of additional cusps. By \cite[Appendix~2]{KS}
complementary series representations satisfy $0<\nu_\tau\leq\theta$
for $\theta=7/64$. Accordingly, the right hand side of \eqref{ineq9}
becomes (cf.\ \eqref{smallu2})
\[\ll_{\eps,a,b,c}
C'_{a,b,c}\,(\ell_1\ell_2)^{1/2+\eps}\,\min(y^{1/2-\theta-\eps},y^{1/2-\theta-b-\eps})\]
for any $0<\eps<1/4-\theta$, where $C'_{a,b,c}$ is a similar
constant as $C_{a,b,c}$ in Theorem~\ref{theorem1} but with a larger
exponent for $\tilde{\lambda}_{\pi_1}+\tilde{\lambda}_{\pi_2}$ and a
sum over larger $d_1+d_2$. For the Eisenstein spectrum the exponent
of $\ell_1\ell_2$ can be lowered to $1/4+\eps$. Similarly, the right
hand side of \eqref{ineq15} becomes, for suitable constants
$A,B,C>0$,
\[\ll_{\eps,k}(\ell_1\ell_2)^A(\tilde{\lambda}_{\pi_1}+\tilde{\lambda}_{\pi_2})^{B+4c}\,|s|^{C+4c},\quad
\frac{1}{2}+\theta+\eps<\Re s<\frac{3}{2}.\] These extensions enable
considerable simplifications in arguments leading to subconvexity
results such as \eqref{burgessbound} which appeared originally in
\cite{BHM} or the more difficult case of Rankin--Selberg
$L$-functions treated in \cite{HM}. In this paper we have decided to
present the theorems in a special case in order to keep the
notational burden minimal and to emphasize the conceptual simplicity
of the approach.
\end{remark}

Selberg \cite{Se} proved in 1965 that for holomorphic cusp forms
(that is, when $\pi_1$ and $\pi_2$ are in the discrete series), the
left hand side of \eqref{spectraldecomposition2} is meromorphic in
$s$ and holomorphic for $\Re(s) > 1/2$ (note that for $\GG$ poles on
the segment $1/2<s\leq 1$ do not occur).\footnote{He adds
regrettingly \cite[p.14]{Se}: ``We cannot make much use of this
function at present [...]"} Good \cite{Go,Go2} was the first to use
the spectral decomposition of shifted convolution sums coming from
the Fourier coefficients of a holomorphic cusp form for estimating
the second moment of the corresponding modular $L$-function on the
critical line. His new key ingredient was to show the polynomial
growth on vertical lines.

The spectral decomposition in the non-holomorphic case has resisted
all attempts so far. Good's method, based on the fact that
holomorphic cusp forms are linear combination of Poincar\'e series,
is not applicable here. Jutila \cite{Ju} and Sarnak \cite{Sa1, Sa2}
independently considered an approximation of the Dirichlet series
\eqref{spectraldecomposition2}, along with a spectral decomposition,
that could be continued to the half plane $\Re s> 1/2$ with a bound
$O(h^{\frac{1}{2}-\sigma+\theta+\eps}|s|^A)$ on vertical lines $\Re
s=\sigma$. However, this approximation introduces an error of
$O(h^{1-\sigma+\eps}|s|^B)$ which one would like to
remove.\footnote{For a striking recent application of the
approximate spectral decomposition see \cite{LLY}. For a careful
analysis of the error term see \cite{Ju3}.} The second author found
in his thesis \cite[Chapter~5]{Ha} that the error signifies missing
harmonics and anticipated that analysis on $\GG\backslash G$ will be
key to obtaining a complete system. Independently, Motohashi
\cite{Mo4,Mo} brought the representation theory of $\GG\backslash G$
and the Kirillov model into the discussion. In this paper, we pursue
his approach further. One of the new ideas we employ is the use of
Sobolev norms for smooth vectors inspired by the recent work of
Venkatesh \cite{Ve} which allows for a soft treatment and avoids the
difficulties imposed by Poincar\'e series and the estimation of
triple product periods. Sobolev norms also play important roles in
related works of Bernstein and Reznikov \cite{BR1,BR,BR3,BR4} and
Cogdell and Piatetski-Shapiro \cite{CoPS}. Shortly after this paper
was finished, Motohashi \cite{Mo3} gave an alternative proof of
Theorem~\ref{theorem2} which makes regularity properties of the
weight functions $F_{k,\tau}(s)$ more transparent.

Finally we remark that there are other important techniques for
understanding shifted convolution sums. The archetype of shifted
convolution sums are the additive divisor sums which have been
studied extensively.\footnote{In fact one can trace back history to
various elegant identities published by Jacobi in 1829.} These
special sums arise from Eisenstein series rather than cusp forms and
their spectral decomposition is very explicitly known by the work of
Motohashi \cite{Mo2} (cf.\ \cite[Lemma~4]{JM}). In the general case
variants of the circle method with the Kloosterman refinement have
been particularly successful \cite{DFI1,DFI2,Ju2,Ha2,HM,Bl,BHM}.
Recently Venkatesh \cite{Ve} developed a geometric method, based on
equidistribution and mixing, which can be applied for shifted
convolution sums of higher rank. In our context, the left hand side
of \eqref{spectraldecomposition} is an automorphic period over a
closed horocycle, so that \cite[Theorem~3.2]{Ve} implies a weaker
but nontrivial version of \eqref{shiftedbound}.

\section{Bounds for the discrete spectrum}

\subsection{Kirillov model and Sobolev norms}

We shall use the notation
\[n(x):=\begin{pmatrix}1&x\\0&1\end{pmatrix},\quad
a(u):=\begin{pmatrix}u&0\\0&1\end{pmatrix},\quad
k(\theta):=\begin{pmatrix}\ \ \ \cos\theta&\sin\theta\\-\sin\theta&
\cos\theta\end{pmatrix},\] and we shall think of these matrices as
elements of $G=\PGL_2(\RR)$. In addition, we shall write
\[e(x):=e^{2\pi ix},\qquad x\in\RR.\]

Let $(\pi,V_\pi)$ be a representation generated by a cusp form on
$\SL_2(\ZZ)\backslash\HH$. Then $(\pi,V_\pi)$ is contained in
$L^2(\GG\backslash G)$, hence it is equipped with a canonical inner
product given by
\begin{equation}\label{inner}\langle\phi_1,\phi_2\rangle:=\int_{\GG\backslash
G}\phi_1(g)\overline{\phi_2(g)}\,dg,\end{equation} where $dg:=dx
(du/u^2) (d\theta/\pi)$ for $g = n(x)a(u)k(\theta)$ is the Haar
measure on $G$. The Kirillov model $\KK(\pi)$ realizes $\pi$ in
$V_{\KK(\pi)}:=L^2(\RR^\times,d^\times u)$ which is equipped with
its own canonical inner product given by
\begin{equation}\label{Winner}\langle W_1,W_2\rangle:=\int_{\RR^\times}
W_1(u)\overline{W_2(u)}\,d^\times u.\end{equation} Let $\|\cdot\|$
denote the norms determined by these inner products.

For a smooth vector $\phi\in V_\pi^\infty$ we define the
corresponding smooth vector $W_\phi\in\ V_{\KK(\pi)}^\infty$ as
\begin{equation}\label{Kirillov}W_\phi(u):=\int_0^1 \phi(n(x)a(u))\,e(-x)\,dx,\qquad
u\in\RR^\times.\end{equation} It gives rise to the Fourier
decomposition
\begin{equation}\label{Fourier}
\phi(n(x)a(u))=\sum_{\substack{n\in\ZZ\\ n\neq
0}}\frac{\lambda_\pi(|n|)}{\sqrt{|n|}}W_\phi(nu)\,e(nx),\qquad
x\in\RR,\quad u \in \RR^{\times},\end{equation} where
$\lambda_\pi(n)$ denotes the $n$-th Hecke eigenvalue of the cusp
form on $\SL_2(\ZZ)\backslash\HH$ that generates $(\pi,V_\pi)$. This
follows from Shalika's multiplicity one theorem combined with
standard facts about Hecke operators, see \cite[Sections~4.1 and
6.2]{CoPS} and \cite[(6.14)--(6.15)]{DFI}. We have the uniform bound
\cite[Proposition~19.6]{DFI}
\begin{equation}\label{DFIbound}
  \sum_{1\leq n \leq x} |\lambda_{\pi}(n)|^2 \ll_{\eps} x (x|\nu_{\pi}|)^{\eps}.
\end{equation}
By Kirillov's theorem \cite[Sections~4.2--4.4]{CoPS}, the canonical
inner products of $V_\pi$ and $V_{\KK(\pi)}$ are related by a
proportionality constant depending only on $\pi$,
\begin{equation}\label{constant}\langle\phi_1,\phi_2\rangle=
C_\pi\langle W_{\phi_1},W_{\phi_2}\rangle,\qquad \phi_1,\phi_2\in
V_\pi^\infty.
\end{equation}
The relations \eqref{Fourier} and \eqref{constant} can also be
verified by classical means, see \cite[Section~4]{DFI} and
\cite[Sections~2 and 4]{BrMo}.

We can evaluate the proportionality constant $C_\pi$ as follows. Let
\begin{displaymath}
  E(g,s) := \frac{1}{2}\sum_{\gamma \in \GG_{\infty}\backslash \GG} U_s(\gamma g)
\end{displaymath}
with $U_s(n(x)a(u)k(\theta)):=|u|^s$ and $\GG_{\infty}:=\{n(x):x \in
\ZZ\}$ denote the standard weight $0$ Eisenstein series on $G$ and
let $\phi\in V_\pi^\infty$ be any vector of pure weight (cf.\
Section~\ref{weightsection}). Then for $\Re s>1$ we have, by the
Rankin--Selberg unfolding technique, \eqref{Fourier}, and Parseval,
\begin{align*}
\langle|\phi(g)|^2, E(g,s)\rangle
&=\frac{1}{2}\int_{\RR^\times}\int_0^1|\phi(n(x)a(u)|^2 \,|u|^{s-2}\,dx\,du\\
&=\frac{1}{2}\sum_{\substack{n\in\ZZ\\ n\neq
0}}\frac{|\lambda_\pi(|n|)|^2}{|n|}
\int_{\RR^\times} |W_\phi(nu)|^2\,|u|^{s-1}\,d^\times u\\
&=\sum_{n=1}^\infty\frac{|\lambda_\pi(n)|^2}{n^s}\int_{\RR^\times}|W_\phi(u)|^2\,|u|^{s-1}\,d^\times
u.
\end{align*}
Taking residues of both sides at $s=1$ yields
\[\|\phi\|^2\res_{s=1}E(g,s)=\|W_{\phi}\|^2\res_{s=1}\frac{L(s,\pi\otimes\tilde\pi)}{\zeta(2s)},\]
Using \eqref{constant} we conclude that\footnote{Here the adjoint square is the same as the symmetric square.}
\[C_\pi=\frac{\vol(\GG\backslash G)}{\zeta(2)}\,L(1,\Ad^2\pi).\] In particular, by
\eqref{DFIbound} and \cite[Theorem~0.2]{HL} we know that
\begin{equation}\label{constantbounds}
\tilde{\lambda}_\pi^{-\eps}\ll_\eps
C_\pi\ll_\eps\tilde{\lambda}_\pi^\eps\end{equation}
for any $\eps > 0$. \\

We shall introduce Sobolev norms for vectors in $V_\pi^\infty$ and
$V_{\KK(\pi)}^\infty$ in terms of the derived action of the Lie
algebra $\fg$ of $G$. We consider the usual basis of $\fg$
consisting of
\[H:=\begin{pmatrix}1&0\\0&-1\end{pmatrix},\quad
R:=\begin{pmatrix}0&1\\0&0\end{pmatrix},\quad
L:=\begin{pmatrix}0&0\\1&0\end{pmatrix},\] and note the commutation
relations $[H,R]=2R$, $[H,L]=-2L$, $[R,L]=H$. Then for smooth
vectors $\phi\in V_\pi^\infty$ and $W\in V_{\KK(\pi)}^\infty$ we
define
\begin{equation}\label{Sobolev}\|\phi\|_{S^{d}}:=\sum_{\ord(\DD)\leq
d}\|\DD\phi\|\qquad\text{and}\qquad
\|W\|_{S^{d}}:=\sum_{\ord(\DD)\leq d}\|\DD W\|,\end{equation} where
$\DD$ ranges over all monomials in $H,R,L$ of order at most $d$ in
the universal enveloping algebra $U(\fg)$. We can see how $U(\fg)$
acts on $\KK(\pi)$: $H$ acts by $2u\frac{d}{du}$ and $R$ acts by
$2\pi i u$, cf.\ \cite[p.155]{Bu}. The Casimir element
$H^2+2RL+2LR=H^2-2H+4RL$ acts by $-4\lambda_\pi$, hence $RL$ acts by
$-\lambda_\pi+u^2\frac{d^2}{du^2}$, therefore $L$ acts by $(2\pi
i)^{-1}(-\lambda_\pi u^{-1}+u\frac{d^2}{du^2})$. This way we obtain
the following estimate \cite[Lemma~8.4]{Ve}:
\begin{equation}\label{akshay}\|W\|_{S^{d}}\ll_d \tilde\lambda_\pi^d\ \|W\|_{A^{2d}}.\end{equation}
Here the norm $\|\cdot\|_{A^{d}}$ was defined in \eqref{akshaynorm}.

\subsection{Normalized Whittaker functions}\label{weightsection} In
order to understand the behavior of $W_\phi(u)$ near $0$ we
decompose $\phi\in V_\pi$ into pure weight pieces
\begin{equation}\label{weights}\phi=\sum_{p\in\ZZ}\phi_p,\end{equation}
where $\phi_p\in V_\pi^\infty$ satisfies
\begin{equation}\label{phip}\phi_p(g
k(\theta))=e^{2ip\theta}\phi_p(g),\qquad g\in
G,\quad\theta\in\RR.\end{equation} Convergence of \eqref{weights} is
understood in $L^2$-norm. Correspondingly, $W_\phi$ decomposes in
$V_{\KK(\pi)}$ as
\begin{equation}\label{Wweights}\
W_\phi=\sum_{p\in\ZZ}W_{\phi_p}.\end{equation} Note that Parseval
gives
\begin{equation}\label{parseval}
\|\phi\|^2=\sum_{p\in\ZZ}\|\phi_p\|^2.
\end{equation}
It is known (cf.\ \cite[(2.14)--(2.27)]{BrMo}) that each
$W_{\phi_p}(u)$ is a constant multiple of some normalized Whittaker
function
\begin{equation}\label{Whittaker}\tilde W_{p,\pi}(u):=\frac{\eps_{p,\pi}(\sgn(u))\, W_{\sgn(u)p,\nu_\pi}(|u|)}
{\left|\GG\left(\frac{1}{2}-\nu_\pi+\sgn(u)p\right)\GG\left(\frac{1}{2}+\nu_\pi+\sgn(u)p\right)\right|^{1/2}},\qquad
u\in\RR^\times,\end{equation} where $\eps_{p,\pi}:\{\pm 1\}\to\{z:
|z|=1\}$ is a suitable phase factor, $W_{\alpha, \beta}$ is the
standard Whittaker function (see \cite[Chapter~XVI]{WW}), and the
right hand side is understood as zero if one of $\frac{1}{2}\pm
\nu_\pi+\sgn(u)p$ is a nonpositive integer. By
\cite[Section~4]{BrMo}, the functions $\tilde
W_{p,\pi}:\RR^\times\to\CC$ form an orthonormal basis of
$L^2(\RR^\times,d^\times u)$, therefore
\begin{equation}\label{scaling}|W_{\phi_p}(u)|=\|W_{\phi_p}\|\,|\tilde W_{p,\pi}(u)|,\qquad
u\in\RR^\times.\end{equation}

We can choose the parameter $\nu_\pi$ so that $\Re\nu_\pi\geq 0$.
Then for $\tilde W_{p,\pi}(u)\neq 0$ we have
\[\frac{\GG\left(\frac{1}{2}+\nu_\pi+\sgn(u)p\right)}{\GG\left(\frac{1}{2}-\nu_\pi+\sgn(u)p\right)}
\ll (|p|+|\nu_\pi|+1)^{2\Re\nu_\pi},\] so that the uniform bound
\cite[(4.5)]{BrMo} (whose proof applies in all cases) yields
\begin{equation}\label{largeu}
\tilde W_{p,\pi}(u)\ll
|u|^{1/2}\left(\frac{|u|}{|p|+|\nu_\pi|+1}\right)^{-1-\Re\nu_\pi}\exp\left(-\frac{|u|}{|p|+|\nu_\pi|+1}\right)
,\quad u\in\RR^\times.\end{equation} If $(\pi,V_\pi)$ belongs to the
principal series (i.e.\ $\Re\nu_\pi=0$) then for $0<\eps<1$ we also
have the uniform bound (cf.\ \cite[(4.3)]{BrMo})
\begin{equation}\label{smallu}
\tilde W_{p,\pi}(u)\ll_\eps (|p|+|\nu_\pi|+1)\,|u|^{1/2-\eps},\qquad
u\in\RR^\times.
\end{equation}
Indeed, for $|u|\geq 1$ the bound \eqref{largeu} is stronger, while
for $|u|<1$ it is an immediate consequence of \cite[(4.2)]{BrMo} and
\cite[Appendix]{HM}. We shall show below that for $0<\eps<1/4$ this
bound holds true even when $(\pi,V_\pi)$ belongs to the discrete
series. We note that representations belonging to the complementary
series (i.e.\ $0<\nu_\pi<1/2$) do not occur in \eqref{spectraldec},
but for completeness we record an analogue of \eqref{smallu} for
this case, valid for $0<\eps<1$ (cf.\ \eqref{largeu},
\cite[(4.2)]{BrMo}, \cite[Appendix]{HM}):
\begin{equation}\label{smallu2}\tilde W_{p,\pi}(u)\ll_\eps (|p|+|\nu_\pi|+1)^{1+\nu_\pi}\,|u|^{1/2-\nu_\pi-\eps},\qquad
u\in\RR^\times.\end{equation}

If $(\pi,V_\pi)$ belongs to the discrete series then
$\nu_\pi=\ell-\frac{1}{2}$, where $\ell\geq 1$ is an integer. For
$|p|<\ell$ we have $\tilde W_{p,\pi}=0$. For $|p|\geq\ell$ it
follows from \eqref{largeu} via $\exp(-t)\ll t^{-1}$ (or from
\cite[(2.16)]{BrMo} by two integration by parts) that
\begin{equation}\label{largeud}
\tilde W_{p,\pi}(u)\ll |3p|^{\ell+3/2}\,|u|^{-\ell-1},\qquad
u\in\RR^\times.\end{equation} The Mellin transform of $W_{p,\pi}(u)$
satisfies the Jacquet--Langlands local functional equation (see
\cite[(4.11)]{BrMo}) which is reflected in the convolution identity
(see \cite[(4.9)]{BrMo})
\[\tilde W_{p,\pi}(u)=(-1)^p\int_0^\infty j_{\ell-\frac{1}{2}}(y)
\,\tilde W_{p,\pi}(y/u)\,d^\times y,\qquad u\in\RR^\times,\] where
\[j_{\ell-\frac{1}{2}}(y):=(-1)^\ell \,2\pi \sqrt{y}\, J_{2\ell-1}(4\pi \sqrt{y}),\qquad y>0.\]
Using the bound \cite[(2.46)]{BrMo} for the kernel
$j_{\ell-\frac{1}{2}}(y)$, we can conclude
\[\tilde W_{p,\pi}(u)\ll \int_0^\infty
\min(y^{1/4},y^\ell)\ |\tilde W_{p,\pi}(y/u)|\,d^\times y,\qquad
u\in\RR^\times.\] We split the integral at $y=|3pu|$ and estimate
the two pieces separately. On the one hand, by Cauchy--Schwarz and
$\|\tilde W_{p,\pi}\|=1$,
\begin{align*}\int_0^{|3pu|}\dots&\leq
\left\{\int_0^{|3pu|}\min(y^{1/2},y^{2\ell})\ d^\times
y\right\}^{1/2} \\&\ll \min(|3pu|^{1/4},|3pu|^\ell).\end{align*} On
the other hand, by \eqref{largeud},
\begin{align*}\int_{|3pu|}^\infty\dots&\ll
|3p|^{\ell+3/2}\int_{|3pu|}^\infty
\min(y^{1/4},y^\ell)\,(y/|u|)^{-\ell-1}\,d^\times y\\ &\ll
|p|^{1/2}\min(|3pu|^{1/4},|3pu|^\ell).\end{align*} All in all we see
that for any $0<\eps<1/4$ we have
\[\tilde W_{p,\pi}(u)\ll |p|^{1/2}\min(|3pu|^{1/4},|3pu|^\ell) \leq |p|^{1/2} |3pu|^{1/2-\eps},\qquad
u\in\RR^\times.\] This implies \eqref{smallu} as claimed.

\subsection{Bounds for smooth vectors}
We can derive a bound for $\|\phi\|_\infty$ in terms of a suitable
Sobolev norm of $\phi$ (cf.\ \eqref{Sobolev}). Let
$y_0:=|p|+|\nu_\pi|+1$. By \eqref{Fourier}, \eqref{phip},
\eqref{scaling}, \eqref{largeu}, and \eqref{smallu} we have
\begin{multline*}
  \|\phi_p\|_\infty \leq\sup_{y>1/2}\sum_{\substack{n\in\ZZ\\ n\neq
0}}\frac{|\lambda_\pi(|n|)|}{\sqrt{|n|}}|W_{\phi_p}(ny)|\\
  \ll_\eps \|W_{\phi_p}\|\sup_{y>1/2} \left\{\sum_{1\leq n \leq y_0/y} \frac{|\lambda_\pi(n)|}{\sqrt{n}} y_0 (ny)^{1/2-\eps} + \sum_{n > y_0/y} \frac{|\lambda_\pi(n)|}{\sqrt{n}} \frac{y_0}{\sqrt{ny}} \,e^{-\frac{ny}{y_0}}\right\}.
\end{multline*}
Together with \eqref{DFIbound}, \eqref{constant}, and
\eqref{constantbounds} we find
\begin{equation*}
\|\phi_p\|_\infty\ll_\eps \|W_{\phi_p} \| \sup_{y > 1/2}
\frac{y_0^2}{\sqrt{y}}(yy_0(1+|\nu_{\pi})|)^{\eps} \ll
(|p|+|\nu_\pi|+1)^{2+\eps}\,\|\phi_p\|.
\end{equation*}
Since $\phi\in V_\pi^\infty$ the norm $\|\phi_p\|$ decays fast in
$p$ which enables to derive an analogue of this inequality for
$\phi$. Indeed, $R-L\in\fg$ acts by the differential operator
$d/d\theta$ (cf.\ \cite[p.155]{Bu}), hence it follows via
\eqref{phip} that
\begin{equation*}
\|\phi_p\|_\infty\ll_\eps
\tilde{\lambda}_\pi^{1+\eps}(1+|p|)^{-1+\eps}\,\|(1+R-L)^3\phi_p\|.\end{equation*}
We apply $(1+R-L)^3\in U(\fg)$ on the weight decomposition
\eqref{weights}, then by Parseval (cf.\ \eqref{parseval}),
Cauchy--Schwarz, and \eqref{weights} we obtain the uniform
bound\footnote{Less explicit versions of this bound were derived by
Bernstein--Reznikov and Venkatesh in more general contexts, without
recourse to Whittaker functions, see \cite[Proposition~4.1]{BR} and
\cite[Lemma~9.3]{Ve}.}
\begin{equation}\label{ineq4}
\|\phi\|_\infty\ll_\eps \tilde{\lambda}_\pi^{1+\eps}\,\|\phi\|_{S^3}.
\end{equation}

In a similar fashion we can derive pointwise bounds for $W_\phi$.
First we combine \eqref{scaling}, \eqref{smallu}, \eqref{constant},
and \eqref{constantbounds} to see that for any $0<\eps<1/4$
\begin{equation*}
W_{\phi_p}(u)\ll_\eps (|p|+|\nu_\pi|+1)^{1+\eps}\,\|\phi_p\| \,
|u|^{1/2-\eps},\qquad u\in\RR^\times,\end{equation*} then we use the
action of $R-L\in\fg$ to conclude that
\begin{equation*}
W_{\phi_p}(u)\ll_\eps
\tilde{\lambda}_\pi^{1/2+\eps}(1+|p|)^{-1+\eps}\,\|(1+R-L)^2\phi_p\| \,
|u|^{1/2-\eps},\quad u\in\RR^\times.\end{equation*} We apply
$(1+R-L)^2\in U(\fg)$ on the weight decomposition \eqref{weights},
then by Parseval (cf.\ \eqref{parseval}), Cauchy--Schwarz, and
\eqref{Wweights} we obtain the uniform bound
\[W_\phi(u)\ll_\eps
\tilde{\lambda}_\pi^{1/2+\eps}\,\|\phi\|_{S^{2}}\,|u|^{1/2-\eps},\qquad
u\in\RR^\times.\] Finally, by replacing $\phi$ by
$(\pm 1+ H^2+2RL+2LR)^cR^bH^a\phi$ we obtain the more general inequality
(cf.\ comments after \eqref{Sobolev})
\begin{multline}\label{Wbound2}\left(u\frac{d}{du}\right)^a
W_\phi(u)\ll_{\eps,a,b,c}\\
\tilde{\lambda}_\pi^{1/2-c+\eps}\,\|\phi\|_{S^{2+a+b+2c}}\,
\min(|u|^{1/2-\eps},|u|^{1/2-b-\eps}),\qquad
u\in\RR^\times,\end{multline} for any $0<\eps<1/4$ and any integers
$a,b,c\geq 0$.

\section{Bounds for the continuous spectrum}

Let $(\pi,V_\pi)$ be a representation generated by an Eisenstein
series on $\SL_2(\ZZ)\backslash\HH$. As $(\pi,V_\pi)$ is not
contained in $L^2(\GG\backslash G)$, we cannot use the definition
\eqref{inner} as in the cuspidal case. Nevertheless, we can still
define the Kirillov model $\KK(\pi)$ and use the definitions
\eqref{Winner}, \eqref{Kirillov}. It turns out that for the purpose
of spectral decomposition the right analogue of \eqref{inner}
reads\footnote{We apologize to the reader that $\pi$ denotes a
constant and a representation at the same time.}
\begin{equation}\label{inner2}\langle\phi_1,\phi_2\rangle:=
\pi^{-1}|\zeta(1+2\nu_\pi)|^2\,\langle
W_{\phi_1},W_{\phi_2}\rangle,\end{equation} then we have
\eqref{constant} and {\it the lower bound part of}
\eqref{constantbounds} with
\[C_\pi:=\pi^{-1}|\zeta(1+2\nu_\pi)|^2.\] Note that $\nu_\pi=0$ does not
occur for a nonzero Eisenstein series. Let $\|\cdot\|$ denote the
norms on $V_\pi$ and $V_{\KK(\pi)}$ determined by these inner
products, then \eqref{Sobolev} defines the corresponding Sobolev
norms $\|\cdot\|_{S^d}$ on $V_\pi^\infty$ and $V_{\KK(\pi)}^\infty$.

For a smooth vector $\phi\in V_\pi^\infty$ we have the following
analogue of \eqref{Fourier}:
\begin{equation}\label{Fourier2}
\phi(n(x)a(u))=W_{\phi,0}(y)+\sum_{\substack{n\in\ZZ\\ n\neq
0}}\frac{\lambda_\pi(|n|)}{\sqrt{|n|}}W_\phi(nu)\,e(nx),\qquad
x\in\RR,\quad u \in \RR^{\times},\end{equation} where
\[W_{\phi,0}(u):=\int_0^1
\phi(n(x)a(u))\,dx,\qquad u\in\RR^\times,\] and
\[\lambda_\pi(n):=\sum_{ab=n}(a/b)^{\nu_\pi}.\]
The equations \eqref{weights}, \eqref{phip}, \eqref{Wweights},
\eqref{parseval} hold true as in the cuspidal case.

It is known (cf.\ \cite[(3.31), (2.16)]{BrMo}) that each
 $W_{\phi_p}(u)$ is a constant multiple of some normalized Whittaker
function of the form \eqref{Whittaker}, therefore we can conclude,
exactly as in the cuspidal case,
\begin{multline}\label{Wbound3}\left(u\frac{d}{du}\right)^a
W_\phi(u)\ll_{\eps,a,b,c}\\
\tilde{\lambda}_\pi^{1/2-c+\eps}\,\|\phi\|_{S^{2+a+b+2c}}\,\min(|u|^{1/2-\eps},|u|^{1/2-b-\eps}),\qquad
u\in\RR^\times,\end{multline} for any $0<\eps<1/4$ and any integers
$a,b,c\geq 0$.

\section{Proof of Theorem~1}

Let $(\pi_1,V_{\pi_1})$ and $(\pi_2,V_{\pi_2})$ be arbitrary
cuspidal automorphic representations of $\GG\backslash G$, and let
$W_{1,2}:\RR^\times\to\CC$ be arbitrary smooth functions of compact
support. There are unique smooth vectors $\phi_i\in
V_{\pi_i}^\infty$ such that $W_{\phi_i}=W_i$. If $h>0$ is an
arbitrary integer and $Y>0$ is arbitrary then we have, by
\eqref{Fourier},
\[\sum_{m+n=h}\frac{\lambda_{\pi_1}(|m|)\lambda_{\pi_2}(|n|)}{\sqrt{|mn|}}
W_1\left(\frac{m}{Y}\right)W_2\left(\frac{n}{Y}\right)=\int_0^1
(\phi_1\phi_2)(n(x)a(Y^{-1}))\,e(-hx)\,dx.\] Let us decompose
spectrally the smooth vector $\phi_1\phi_2\in L^2(\GG\backslash G)$
according to \eqref{spectraldec},
\[\phi_1\phi_2=\int_\tau\psi_\tau\,d\tau.\] This
decomposition is unique and converges in the topology defined by the
Sobolev norms $\|\cdot\|_{S^d}$, see Propositions~1.3 and 1.4 in
\cite{CoPS}. In particular, $\psi_\tau\in V_\tau^\infty$, and the
decomposition is compatible with the actions of $G$ and $\fg$. The
explicit knowledge of the projections $L^2(\GG\backslash G)\to
V_\tau$ yields (cf.\ \cite[(2.14), (3.31), Lemma~2]{BrMo}), in
combination with Plancherel, \eqref{inner}, and \eqref{inner2},
\begin{equation}\label{Plancherel}\|\DD(\phi_1\phi_2)\|^2=\int_\tau \|\DD \psi_\tau\|^2\,d\tau,\qquad
\DD\in U(\fg).\end{equation} By Corollary to Lemma~1.1 in
\cite{CoPS}, the functional $\phi\mapsto W_\phi(h)$ is continuous in
the topology defined by the Sobolev norms $\|\cdot\|_{S^d}$ (this
also follows from \eqref{Wbound2}, \eqref{Wbound3}, and Plancherel),
hence the above imply, in combination with \eqref{Fourier} and
\eqref{Fourier2},
\[\sum_{m+n=h}\frac{\lambda_{\pi_1}(|m|)\lambda_{\pi_2}(|n|)}{\sqrt{|mn|}}
W_1\left(\frac{m}{Y}\right)W_2\left(\frac{n}{Y}\right)=
\int_{\tau\neq\tau_0}\frac{\lambda_\tau(h)}{\sqrt{h}}W_{\tau}\left(\frac{h}{Y}\right)d\tau\]
where $W_\tau:=W_{\psi_\tau}$. This is just
\eqref{spectraldecomposition}. On the right hand side we have, by
\eqref{Wbound2} and \eqref{Wbound3},
\begin{multline*}\tilde{\lambda}_\tau^{c}\left(y\frac{d}{dy}\right)^{a}W_{\tau}(y)\ll_{\eps,a,b,c}\\
\tilde{\lambda}_\tau^{-3/2+\eps}\|\psi_\tau\|_{S^{6+a+b+2c}}\,\min(y^{1/2-\eps},y^{1/2-b-\eps}),\qquad
y>0.\end{multline*} A combination of Cauchy--Schwarz, Weyl's law,
and \eqref{Plancherel} then shows
\begin{multline}\label{ineq8}\int_{\tau\neq\tau_0}
\tilde{\lambda}_\tau^{c}\left|\left(y\frac{d}{dy}\right)^{a}W_{\tau}(y)\right|\,d\tau
\ll_{\eps,a,b,c}\\
\|\phi_1\phi_2\|_{S^{6+a+b+2c}}\,\min(y^{1/2-\eps},y^{1/2-b-\eps}),\qquad
y>0.\end{multline} Using \eqref{ineq4} and the Leibniz rule for
derivations we see that
\[\|\phi_1\phi_2\|_{S^{d}}\ll_{\eps,d}
(\tilde{\lambda}_{\pi_1}+\tilde{\lambda}_{\pi_2})^{1+\eps}\sum_{d_1+d_2=d+3}\|\phi_1\|_{S^{d_1}}\,\|\phi_2\|_{S^{d_2}}\]
for any integer $d \geq 0$, which implies, by \eqref{constant},
\eqref{constantbounds}, and \eqref{akshay},
\[\|\phi_1\phi_2\|_{S^{d}}\ll_{\eps,d}
(\tilde{\lambda}_{\pi_1}+\tilde{\lambda}_{\pi_2})^{d+4+\eps}\sum_{d_1+d_2=2d+6}\|W_1\|_{A^{d_1}}\,\|W_2\|_{A^{d_2}}.\]
We combine this inequality with \eqref{ineq8} to arrive at
\eqref{ineq9}.

In retrospect we can see that this proof works for all test
functions $W_{1,2}:\RR^\times\to\CC$ whose norms
$\|W_{1,2}\|_{A^{d}}$ exist for $d=18+2a+2b+4c$.

\section{Proof of Theorem~2}

Let $(\pi_1,V_{\pi_1})$ and $(\pi_2,V_{\pi_2})$ be arbitrary
cuspidal automorphic representations of $\GG\backslash G$, and let
$c,k\geq 0$ be arbitrary integers. Let $q>0$ be an integer to be
determined later in terms of $c$ and $k$. Consider the function
$G:[0,\infty)\to\CC$ defined by
\[G(t):=\begin{cases}\{t(1-t)\}^q,& 0\leq t\leq 1,\\
0,&t>1,\end{cases}\] and its Laplace transform (composed with
$z\mapsto -z$) defined by
\[\tilde G(z):=\int_0^\infty G(t)\, e^{zt}\,dt,\qquad z\in\CC.\]
Note that $\tilde G(z)$ is entire and by successive integration by
parts it satisfies the uniform bound
\begin{equation}\label{ineq11}\tilde G(z)\ll_q |z|^{-q-1},\qquad
\Re z=1.\end{equation}

If $m,n\geq 1$ are arbitrary integers and $Y>0$ is arbitrary then by
\eqref{ineq11} and the theory of the Laplace transform we have the
identity
\begin{equation}\label{eq1}\left(\frac{mn}{Y^2}\right)^{k/2}G\left(\frac{m+n}{Y}\right)=
\frac{1}{2\pi i}\int_{(1)}\tilde G(z)\,
W_k\left(\frac{m}{Y},z\right)W_k\left(\frac{n}{Y},z\right)\,dz,\end{equation}
where $W_k:[0,\infty)\times\CC\to\CC$ is the function defined by
\begin{equation}\label{weightz}W_k(t,z):=t^{k/2}e^{-zt},\qquad t\geq 0,\quad
z\in\CC.\end{equation} By Theorem~\ref{theorem1} we can see that if
$k$ is sufficiently large in terms of $c$ then there exist functions
$W_{k,\tau}:(0,\infty)\times\{z:\Re z>0\}\to\CC$ depending only on
$\pi_{1,2}$, $k$, and $\tau$ such that the following two properties
hold for all $z$ with $\Re z>0$. If $h>0$ is an arbitrary integer
and $Y>0$ is arbitrary then we have the decomposition over the full
spectrum (excluding the trivial representation)
\begin{equation}\label{spectralz}\sum_{\substack{m,n\geq 1\\ m-n=h}}
\frac{\lambda_{\pi_1}(m)\lambda_{\pi_2}(n)}{\sqrt{mn}}
W_k\left(\frac{m}{Y},z\right)W_k\left(\frac{n}{Y},z\right)=
\int_{\tau\neq\tau_0}\frac{\lambda_\tau(h)}{\sqrt{h}}W_{k,\tau}\left(\frac{h}{Y},z\right)d\tau,\end{equation}
and we have the uniform bounds
\[\int_{\tau\neq\tau_0}\tilde{\lambda}_\tau^{c}
\left|W_{k,\tau}(y,z)\right|\,d\tau \ll_{\eps,c}
C_{0,1,c}(z)\,\min(y^{1/2-\eps},y^{-1/2-\eps}),\quad y>0,\quad\Re
z>0,\] where $0<\eps<1/4$ is arbitrary and
\[C_{0,1,c}(z):=(\tilde{\lambda}_{\pi_1}+\tilde{\lambda}_{\pi_2})^{12+2c}
\sum_{d_1+d_2=20+4c}\|W_k(\cdot,z)\|_{A^{d_1}}\,\|W_k(\cdot,z)\|_{A^{d_2}}.\]
Here we use the convention that $W_k(u,z)=0$ for $u<0$. By
\eqref{akshaynorm} and \eqref{weightz} the right hand side exists as
long as $k>60+12c$. Under this condition \eqref{spectralz} is
justified and we conclude that
\begin{multline}\label{ineq20}\int_{\tau\neq\tau_0}\tilde{\lambda}_\tau^{c}
\left|W_{k,\tau}(y,z)\right|\,d\tau \ll_{\eps,k}\\
(\tilde{\lambda}_{\pi_1}+\tilde{\lambda}_{\pi_2})^{12+4c}\,
|z|^{20+4c}\,\min(y^{1/2-\eps},y^{-1/2-\eps}),\qquad y>0,\quad\Re
z=1.\end{multline}

The functions $W_{k,\tau}(y,z)$ furnished by the proof of
Theorem~\ref{theorem1} also have good behavior for individual
$\tau$. To see this, denote by $\phi_{i,z}\in V_{\pi_i}$ the two
vectors corresponding to the left hand side of \eqref{spectralz} in
the proof of Theorem~\ref{theorem1}. By \eqref{Wbound2} and
\eqref{Wbound3} we have
\[W_{k,\tau}(y,z)\ll_{\eps,\tau}\|\psi_{\tau,z}\|_{S^3}\,\min(y^{1/2-\eps},y^{-1/2-\eps}),\qquad y>0,\quad\Re(z)>0,\]
where $\psi_{\tau,z}\in V_\tau$ is the projection of
$\phi_{1,z}\phi_{2,z}$ on $V_\tau$. Applying \eqref{constant} and
\eqref{inner2}, then \eqref{akshay},
\[\|\psi_{\tau,z}\|_{S^3}\ll_\tau \|W_{k,\tau}(\cdot,z)\|_{S^3}\ll_\tau
\|W_{k,\tau}(\cdot,z)\|_{A^{6}},\qquad\Re(z)>0,\] whence by
\eqref{weightz} we can conclude
\begin{equation}\label{Wkpointwise}W_{k,\tau}(y,z)\ll_{\eps,k,\tau}|z|^{6}\,\min(y^{1/2-\eps},y^{-1/2-\eps}),\qquad
y>0,\quad\Re z=1.\end{equation} In an almost identical fashion,
\begin{equation}\label{Wkpointwise2}y\frac{d}{dy}W_{k,\tau}(y,z)\ll_{\eps,k,\tau}|z|^{6}\,y^{1/2-\eps},\qquad
y>0,\quad\Re z=1.\end{equation} Finally, by \eqref{Wbound2},
\eqref{Wbound3}, \eqref{constant}, \eqref{inner2}, and
\eqref{akshay},
\begin{align*}W_{k,\tau}(y,z)-W_{k,\tau}(y,z')&\ll_{\tau,y}
\|\psi_{\tau,z}-\psi_{\tau,z'}\|_{S^{2}}\\
&\ll_{\tau,y}
\|W_{k,\tau}(\cdot,z)-W_{k,\tau}(\cdot,z')\|_{S^{2}}\\
&\ll_{\tau,y}
\|W_{k,\tau}(\cdot,z)-W_{k,\tau}(\cdot,z')\|_{A^{4}},\qquad \Re z,\
\Re z'>0,\end{align*} whence by \eqref{weightz} we can conclude that
\begin{equation}\label{Wkcont}\lim_{z'\to z}W_{k,\tau}(y,z')=W_{k,\tau}(y,z),\qquad y>0,\quad \Re
z>0.\end{equation}

By \eqref{ineq11}, \eqref{Wkpointwise}, \eqref{Wkpointwise2}, and
\eqref{Wkcont}, the integral
\begin{equation}\label{Hdef}H_{k,\tau}(y):=\frac{1}{2\pi i}\int_{(1)}\tilde G(z) \,W_{k,\tau}(y,z)\,dz,\qquad
y>0,\end{equation} defines a differentiable function
$H_{k,\tau}:(0,\infty)\to\CC$ for all $\tau$ and satisfies the
uniform bound
\begin{equation}\label{Hpointwise}
H_{k,\tau}(y)\ll_{\eps,k,\tau}\min(y^{1/2-\eps},y^{-1/2-\eps}),\qquad
y>0,\end{equation} as long as $q>6$. Then \eqref{ineq11},
\eqref{eq1}, \eqref{weightz}, \eqref{spectralz}, \eqref{ineq20}, and
two applications of Fubini's theorem show that
\[\sum_{\substack{m,n\geq 1\\ m-n=h}}
\frac{\lambda_{\pi_1}(m)\lambda_{\pi_2}(n)}{\sqrt{mn}}\left(\frac{mn}{Y^2}\right)^{k/2}G\left(\frac{m+n}{Y}\right)=
\int_{\tau\neq\tau_0}\frac{\lambda_\tau(h)}{\sqrt{h}}H_{k,\tau}\left(\frac{h}{Y}\right)d\tau,\]
as long as $q>20+4c$. We specify $q:=21+4c$ and evaluate the Mellin
transform in $Y$ at $1-s$ of both sides:
\begin{equation}\begin{split}\label{eq5}&\int_0^\infty Y^{1-s-k}\sum_{\substack{m,n\geq
1\\ m-n=h}}
\lambda_{\pi_1}(m)\lambda_{\pi_2}(n)\,(mn)^\frac{k-1}{2}\,G\left(\frac{m+n}{Y}\right)\frac{dY}{Y}\\
=&\int_0^\infty
Y^{1-s}\int_{\tau\neq\tau_0}\frac{\lambda_\tau(h)}{\sqrt{h}}H_{k,\tau}\left(\frac{h}{Y}\right)d\tau\,\frac{dY}{Y}.
\end{split}\end{equation}
The left hand side is absolutely convergent for $\Re s>1$, and by
Fubini's theorem it equals
\[\hat G(s-1)\sum_{\substack{m,n\geq 1\\ m-n=h}}
\frac{\lambda_{\pi_1}(m)\lambda_{\pi_2}(n)\,(mn)^\frac{k-1}{2}}{(m+n)^{s+k-1}},\qquad\Re(s)>1.\]
The right hand side is absolutely convergent for $\frac{1}{2}<\Re
s<\frac{3}{2}$ by \eqref{Hdef}, \eqref{ineq11}, \eqref{ineq20}, and
by Fubini's theorem it equals
\[h^{\frac{1}{2}-s}\int_{\tau\neq\tau_0}\lambda_{\tau}(h)\,\hat H_{k,\tau}(s-1)\,d\tau,\qquad \frac{1}{2}<\Re s<\frac{3}{2}.\]
The Mellin transforms $\hat H_{k,\tau}(s-1)$ are holomorphic
functions in the strip $\frac{1}{2}<\Re s<\frac{3}{2}$ by
\eqref{Hpointwise}, and by \eqref{Hdef}, \eqref{ineq11},
\eqref{ineq20} they satisfy the uniform bound
\begin{equation}\label{HMellin}\int_{\tau\neq\tau_0}\tilde{\lambda}_\tau^{c}\,
\bigl|\hat
H_{k,\tau}(s-1)\bigr|\,d\tau\ll_{\eps,k}(\lambda_{\pi_1}+\lambda_{\pi_2})^{12+4c},\qquad
\frac{1}{2}+\eps<\Re s<\frac{3}{2}.\end{equation} Finally we observe
that
\[\hat G(s-1)=\int_0^1 \{t(1-t)\}^q\, t^{s-2}\,dt=\frac{\GG(s+q-1)\,\GG(q+1)}{\GG(s+2q)},\qquad\Re s>1-q,\]
therefore
\[\hat G(s-1)\gg_q |s|^{-q-1},\qquad \frac{1}{2}<\Re s<\frac{3}{2}.\]
We put
\[F_{k,\tau}(s):=\frac{\hat H_{k,\tau}(s-1)}{\hat G(s-1)},\qquad \frac{1}{2}<\Re s<\frac{3}{2},\]
then the statements of Theorem~\ref{theorem2} are immediate. In
particular, \eqref{spectraldecomposition2} follows from the
comparison of the two sides in \eqref{eq5}, while \eqref{ineq15} is
a consequence of \eqref{HMellin}.

\section{Acknowledgements}

The authors are indebted to Akshay Venkatesh for sharing with them
his ideas and knowledge generously. An important initial discussion
took place at the workshop ``Subconvexity bounds for $L$-functions"
in October 2006, and the authors are grateful to the organizers and
the American Institute of Mathematics for invitation and funding.
The referees provided valuable comments, their efforts are
appreciated as well. Last but not least, Gergely Harcos would like
to thank Peter Sarnak for his constant support, especially for
emphasizing him the importance of the problem solved in the present
paper and for showing him the beauty of automorphic forms.

\end{document}